\documentclass{amsart}
 \usepackage{graphicx}
 \usepackage{amssymb}

\newtheorem{theorem}{Theorem}[section]
\newtheorem{lemma}[theorem]{Lemma}

\theoremstyle{definition}
\newtheorem{definition}[theorem]{Definition}
\newtheorem{example}[theorem]{Example}

\renewcommand{\mod}{{\bf mod \/}}

\theoremstyle{remark}

\numberwithin{equation}{section}

\begin{document}

\title{There are Infinitely Many Pairs of Twin Prime}

\author{Zhanle Du}
\address{Chinese Academy of Sciences, A20 Datun rd.
 Chaoyang Dst. Beijing 100012, China}
%\curraddr{Current }
\email{zldu@bao.ac.cn}
\thanks{All rights reserved. This is one of a serial works, and was
supported financially by Pro. Wei WANG and on loan in part.
\\ Corresponding author: Zhanle Du.}

\author{Shouyu Du}
\address{Chinese Academy of Sciences, 99 Donggang rd.
 ShiJiaZhuang, HeBei, 050031, China}
\email{shouyudu@yahoo.com.cn}

%    General info
\subjclass{11A41}
% amsart changed \newcommand*\subjclass[2][1991]{%--------to \newcommand*\subjclass[2][2000]{%

\date{Oct.3, 2005}%January 1, 1994 and, in revised form, June 22, 1994.}

%\dedicatory{This paper is dedicated to our authors.}

\keywords{prime,twin primes, floor function, ceiling function,
integral operator}%sieve,

\begin{abstract}
We proved that there are infinitely many pairs of twin prime.
\end{abstract}

\maketitle

\section{Introduction}
\label{sect:intro}

Let P=$\{p_1,p_2,...,p_v \}= \{2,3,...,p_v \}$ be the primes not
exceeding $\sqrt{n}$, then the number of primes not exceeding $n$
\cite{Rosen} is,
\begin{equation}
    \label{equ:prime1}
        \pi(n)=\left\{\begin{array}    {lr}
            (\pi(\sqrt{n})-1)+n
            -\left(\left\lfloor \frac{n}{p_1}\right\rfloor+\left\lfloor \frac{n}{p_2}\right\rfloor
            +\cdots +\left\lfloor \frac{n}{p_v}\right\rfloor \right) \\
            +\left(\left\lfloor \frac{n}{p_1p_2}\right\rfloor+\left\lfloor \frac{n}{p_1p_3}\right\rfloor
            +\cdots +\left\lfloor \frac{n}{p_{v-1}p_v}\right\rfloor \right)
            %\\-\left(\left\lfloor \frac{n}{p_1p_2p_3}\right\rfloor+\left\lfloor \frac{n}{p_1p_2p_4}\right\rfloor
            %+\cdots +\left\lfloor \frac{n}{p_{v-2}p_{v-1}p_v}\right\rfloor \right) \\
            +\cdots.
         \end{array}
         \right.
\end{equation}
For simplicity, we can write it as,
\begin{equation}
    \label{equ:prime2}
    \begin{array}{rl}
            \pi(n) &
            =(\pi(\sqrt{n})-1)+n\left[1- \frac{1}{p_1}\right] \left[1- \frac{1}{p_2}\right]
            \cdots \left[1- \frac{1}{p_v}\right]\\
            & =(\pi(\sqrt{n})-1)+n\prod_{i=1}^{v}\left[1- \frac{1}{p_i}\right], \\
         \end{array}
\end{equation}
where  $n\left[1-
\frac{1}{p_i}\right]=n-\left\lfloor\frac{n}{p_i}\right\rfloor$,
$n\left[1- \frac{1}{p_i}\right]\left[1-
\frac{1}{p_j}\right]=n-\left\lfloor\frac{n}{p_i}\right\rfloor
-\left\lfloor\frac{n}{p_j}\right\rfloor+\left\lfloor\frac{n}{p_ip_j}\right\rfloor$.
The operator $\left[1- \frac{1}{p_i}\right]$ will leave the items
which are not multiples of $p_i$.

In this paper, $\left\lfloor x\right\rfloor\leq x$ is the floor
function, $\left\lceil x\right\rceil\geq x$  the ceiling function
of $x$.  The integral operator `$[\ ]$', which is not the floor
function in this paper, has meanings only operating on (real)
number: $m[\ ]=[m]=\lfloor m\rfloor$.

If $p_t$ and $p_t+2$ are both primes, then this prime pair is
called twin primes. Is there infinitely twin primes? It is a
long-term unsolved
conjecture\cite{Rosen,Hardy,Indlekofer96,Suzuki}. We will proved
it as theorem \ref{theorem}.

\begin{theorem}
\label{theorem} There are infinitely many pairs of twin prime.
\end{theorem}

\section{ The number of twin primes from $p_v^2$ to $p_{v+1}^2$}
\label{sect:Formula}

Let $Z=\{1,2,\cdots,n \}(n<p_{v+1}^2)$ be a natural arithmetic
progression, $Z'=Z+2=\{3,5,\cdots,n+2 \}$ be its accompanying
arithmetic progression, so that $Z'_k=Z_k+2, k=1,2,\cdots, n$.
There are $n$ such pairs.
\begin{equation}
    \label{equ:setZ}\left\{
    \begin{array}{llccccrr}
            Z  & =\{ &   1, &   2, & \cdots, &   n  & \}\\
            Z' & =\{ &   3, &   5, & \cdots, &   n+2 & \}.\\
         \end{array} \right.
\end{equation}

If we delete all the pairs in which one or both items $Z_k$ and
$Z'_k$ are composite integers of all primes $p_i\leq \sqrt{n}$,
then the pairs left are all twin primes.

For a given $p_i$, first we delete the multiples of $p_i$ in set
$Z$, or the items of $Z_k \mod p_i=0$,
\begin{equation}
    \label{equ:y1}
    y(p_i)=\left\lfloor\frac{n}{p_i}\right\rfloor \equiv
    n\left[\frac{1}{p_i}\right].
\end{equation}
Secondly we delete the multiples of $p_i$ in $Z'$, i.e. the items
of $Z'_k \mod p_i=0$, or $Z_k \mod p_i= p_i-2$,
\begin{equation}
    \label{equ:y2}
    y'(p_i)=\left\lfloor\frac{n+2}{p_i}\right\rfloor \equiv
    n\left[\frac{\widetilde{1}}{p_i}\right].
\end{equation}
Because $Z_k\mod p_i\ne Z'_k\mod p_i$ when $p_i\not=2$, the items
deleted by $y(p_i)$ and $y'(p_i)$ are not the same. For $p_i=2$,
$Z_k\mod p_i= Z'_k\mod p_i$, we  need only delete the items in set
$Z$:
\begin{equation}
    \label{equ:yfor2}
    y'(p_i)=0, \quad\mbox{for}\quad p_i=2.
\end{equation}

Eq. (\ref{equ:y1}) will delete all pairs with $Z_k \mod p_i=0$ in
set $Z$, and Eq. (\ref{equ:y2}) will delete all pairs with $Z_k
\mod p_i= p_i-2$ in set $Z$.
\begin{equation}
    \label{equ:y3}
    y'(p_i)=
    n\left[\frac{\widetilde{1}}{p_i}\right]=n\left\lfloor\frac{1}{p_i}\right\rfloor+\delta,
\end{equation}
where $0 \leq \delta \leq 1$ with,
\begin{equation}
  \label{equ:delta}
        \delta=\left\{\begin{array}{lr}
            0: \quad 0\leq n \mod p_i+2 < p_i\\
            1: \quad  \mbox{else}.\\ %m {\ \bf mod\ } p_i+(p_i-1)-(2n-1) {\ \bf mod\ } p_i \geq p_i.\\
         \end{array}
         \right.
\end{equation}
The pairs left, when deleted all the multiples of $p_i$ in both
$Z$ and $Z'$, have,
\begin{equation}
\label{equ:M}
    M(p_i)=n-y(p_i)-y'(p_i)\equiv
    n\left[1-\frac{1}{p_i}-\frac{\widetilde{1}}{p_i}\right].
\end{equation}
The operator
$\left[1-\frac{1}{p_i}-\frac{\widetilde{1}}{p_i}\right]$, when
operating on $n$, will leaves the items having not those of $Z_k
\mod p_i=0,p_i-2$.

For another prime $p_j$, we should also delete its multiples in
both $Z$ and $Z'$, but if there is an item which is composite
integer of both $p_i$ and $p_j$ then this item should be deleted
only once. Thus
$n\left[\frac{1}{p_j}\right]\left[1-\frac{1}{p_{i}}-\frac{\widetilde{1}}{p_{i}}\right]$
will delete the multiples of $p_j$ in $Z$, which have not been
deleted by $p_i$. $n\left[\frac{\widetilde{1}}{p_j}\right]
\left[1-\frac{1}{p_{i}}-\frac{\widetilde{1}}{p_{i}}\right]$ will
delete the multiples of $p_j$ in $Z'$, which have not been deleted
by $p_i$. The pairs will leave,
\begin{equation}
  \label{equ:expandp0}
    \begin{array}{rl}
    M(p_i,p_j) & = n\left[1-\frac{1}{p_i}-\frac{\widetilde{1}}{p_i}\right]
                   -n\left[\frac{1}{p_j}\right]\left[1-\frac{1}{p_i}-\frac{\widetilde{1}}{p_i}\right]
                   -n\left[\frac{\widetilde{1}}{p_j}\right]
                   \left[1-\frac{1}{p_i}-\frac{\widetilde{1}}{p_i}\right]\\
               & \equiv n\left[1-\frac{1}{p_i}-\frac{\widetilde{1}}{p_i}\right]
                   \left[1-\frac{1}{p_j}-\frac{\widetilde{1}}{p_j}\right]\\
               & = n\left[1-\frac{1}{p_i}-\frac{\widetilde{1}}{p_i}
                   -\frac{1}{p_j}+\frac{1}{p_ip_j}+\frac{\widetilde{1}}{p_i}\frac{1}{p_j}
                   -\frac{\widetilde{1}}{p_j}+\frac{1}{p_i}\frac{\widetilde{1}}{p_j}
                   +\frac{\widetilde{1}}{p_i}\frac{\widetilde{1}}{p_j}
                   \right],
    \end{array}
\end{equation}
the meaning is as follows,
\begin{equation}
  \label{equ:expandp}
  \left\{
   \begin{array}{lr}
   n\left[\frac{1}{p_i}\right]=\left\lfloor\frac{n}{p_i}\right\rfloor\\
   n\left[\frac{\widetilde{1}}{p_i}\right]=\left\lfloor\frac{n+2}{p_i}\right\rfloor\\
   n\left[\frac{1}{p_i}\right]\left[\frac{1}{p_j}\right]
       =n\left[\frac{1}{p_ip_j}\right]
       %=\left[\frac{\left[\frac{m}{p_i}\right]}{p_j}\right]
       =\left\lfloor\frac{n}{p_ip_j}\right\rfloor\\
   n\left[\frac{\widetilde{1}}{p_i}\right]\left[\frac{\widetilde{1}}{p_j}\right]
       =n\left[\frac{\widetilde{1}}{p_ip_j}\right]
       =\left\lfloor\frac{n+2}{p_ip_j}\right\rfloor\\
   n\left[\frac{1}{p_i}\right]\left[\frac{\widetilde{1}}{p_j}\right]
     =n\left[\frac{1}{p_i}\frac{\widetilde{1}}{p_j}\right]
       %=\left[\frac{\left[\frac{m}{p_i}\right]+\theta_{j}}{p_j}\right]
       =\left\lfloor\frac{n+\theta_{i,j}}{p_ip_j}\right\rfloor
       =\left\lfloor\frac{n+p_ip_j-\lambda_{i,j}}{p_ip_j}\right\rfloor,\\
 \end{array}
 \right.
\end{equation}
where $1\leq\lambda_{i,j}\leq p_ip_j$ is the first (number)
position in $Z$ with $p_i|Z_{\lambda_{i,j}}$ and
$p_j|Z'_{\lambda_{i,j}}$,
\begin{equation}
  \label{equ:lambda}
  \left\{
   \begin{array}{l}
         Z_{\lambda_{i,j}}\mod p_i=0\\
         Z_{\lambda_{i,j}}\mod p_j= p_j-2.\\
 \end{array}
 \right.
\end{equation}
Let $X=\{\lambda p_i,
\lambda=1,2,\cdots,\left\lfloor\frac{n}{p_i}\right\rfloor\}$,
$\lambda_j$ be the first item with $\lambda_j p_i\mod p_j=p_j-2$,
then $\lambda_{i,j}=\lambda_{j}p_i$,
\begin{equation}
  \label{equ:lambda2}
   \begin{array}{l}
     n\left[\frac{1}{p_i}\right]\left[\frac{\widetilde{1}}{p_j}\right]
       =\left\lfloor\frac{n+p_ip_j-\lambda_{i,j}}{p_ip_j}\right\rfloor
       =\left\lfloor\frac{
           \left\lfloor\frac{n}{p_i}\right\rfloor
            +p_j-\lambda_{j}
         }{p_j}\right\rfloor.\\
 \end{array}
\end{equation}
If there is no such $\lambda_{i,j}$ in $Z$, i.e.,
$\lambda_{i,j}\geq (n+1)$, then the last item in Eq.
(\ref{equ:expandp}) will equal zero. Because $1\leq \lambda_j\leq
(p_j-1)$, so $p_i\leq\theta_{i,j}=p_ip_j-\lambda_{i,j}\leq
p_i(p_j-1)$.

When the multiples of all primes $p_i\leq p_v <\sqrt{n}$ in both
$Z$ and $Z'$ have been deleted, the pairs left will be prime pairs
and have,
\begin{equation}
  \label{equ:D0n}
   \begin{array}{rl}
    D_0(n) & =n\left[1-\frac{1}{p_1}\right]
            \left[1-\frac{1}{p_2}-\frac{\widetilde{1}}{p_2}\right]
            \cdots\left[1-\frac{1}{p_{v}}-\frac{\widetilde{1}}{p_{v}}\right]\\
        & = n\left[1-\frac{1}{p_1}\right]
           \prod_{i=2}^v\limits\left[1-\frac{1}{p_i}-\frac{\widetilde{1}}{p_i}\right],\\
 \end{array}
\end{equation}

The total number of twin primes in $n$ is,
\begin{equation}
  \label{equ:Dn}
   \begin{array}{rl}
    D(n) & =D_0(n)+D(\sqrt{n})-D_1\\
         & =n\left[1-\frac{1}{p_1}\right]
            \prod_{i=2}^v\limits\left[1-\frac{1}{p_i}-\frac{\widetilde{1}}{p_i}\right]
            +D(\sqrt{n})-D_1,\\
 \end{array}
\end{equation}
where $D(\sqrt{n})\geq0$ is the number of twin primes
$[p_i,p_{i+1}]$ when $p_i\leq p_v<\sqrt{n}$, where
$p_{i+1}=p_{i}+2$ and,
\begin{equation}
  \label{equ:D1}
        D_1=\left\{\begin{array}{lr}
            1 \quad \mbox{if}\ n+2=p_{v+1}^2\ \mbox{and}\ n\mbox{ }=prime\\
            0 \quad \mbox{else}.\\
         \end{array}
         \right.
\end{equation}
Because $n<p_{v+1}^2$, if $n=p_{v+1}^2-1$ then $n\mod 2=0$; if
$n=p_{v+1}^2-2$ and $n$ is prime then $[n,n+2]$ is not twin prime.

Besides, the pair [1,3] has been already deleted since $1\mod
p_2=p_2-2$ for $p_2=3$.

     %The twin prime [p,p+2] satisfies,
     %\begin{equation}
     %  \label{equ:Dp61}
     %    \left\{
     %        \begin{array}{rl}
     %            p \mod 2     & \ne0\\
     %            p \mod 3 & \ne0\\
     %            (p+2) \mod 3 & \ne0,\\
     %        \end{array}
     %    \right.
     %\end{equation}
     %that is,
     %\begin{equation}
     %  \label{equ:Dp62}
     %    \left\{
     %        \begin{array}{rl}
     %            p \mod 2 & =1\\
     %            p \mod 3 & =2,\\
     %        \end{array}
     %    \right.
     %\end{equation}
     %or
     %\begin{equation}
     %  \label{equ:Dp63}
     %            p \mod 6  =5.\\
     %\end{equation}
     %In other words, the twin prime must be the form [6k-1,6k+1], k is
     %integer.

\begin{example}
$n=41$, $Z=[1,2,\cdots,41]$, $p_i=[2,3,5]$,

$\overrightarrow{Z_k \mod 2\ne0}\ [1,3,\cdots, 41]$

$\overrightarrow{Z_k\mod 3\ne0}\
[1,5,7,11,13,17,19,23,25,29,31,35,37,41]$

$\overrightarrow{Z_k\mod 5\ne0}\ [1,7,11,13,17,19,23,29,31,37,41]$

$\overrightarrow{Z_k \mod 3\ne1}\ [11,17,23,29,41]$

 $\overrightarrow{Z_k \mod 5\ne3}\ [11,17,29,41]$.

 So $D_0(41)=4$.

We can check it from Eq. (\ref{equ:D0n}) directly.
\begin{displaymath}
  \label{equ:exp1}
  \!\!\!\!\!\!
  \begin{array}{rl}
   D_0(41)\!\!\!\!\!\! & =n\left[1-\frac{1}{2}\right]
            \left[1-\frac{1}{3}-\frac{\widetilde{1}}{3}\right]
            \left[1-\frac{1}{5}-\frac{\widetilde{1}}{5}\right]\\
       & =41\left[1-\frac{1}{2}-\frac{1}{3}+\frac{1}{6}
        -\frac{\widetilde{1}}{3}+\frac{1}{2}\frac{\widetilde{1}}{3}
        -\frac{1}{5}+\frac{1}{10}+\frac{1}{15}-\frac{1}{30}
        +\frac{1}{5}\frac{\widetilde{1}}{3}-\frac{1}{10}\frac{\widetilde{1}}{3}
        \right.\\
       & {\ \ }\left.
        -\frac{\widetilde{1}}{5}+\frac{1}{2}\frac{\widetilde{1}}{5}
        +\frac{1}{3}\frac{\widetilde{1}}{5}
        -\frac{1}{6}\frac{\widetilde{1}}{5}
        +\frac{\widetilde{1}}{15}-\frac{1}{2}\frac{\widetilde{1}}{15}\right]\\
       & =41-\left\lfloor\frac{41}{2} \right\rfloor
        -\left\lfloor\frac{41}{3} \right\rfloor+\left\lfloor\frac{41}{6} \right\rfloor
        -\left[\frac{\widetilde{41}}{3} \right]
        +\left[\frac{41}{2}\frac{\widetilde{1}}{3}\right]\\
       & {\ \ }-\left\lfloor\frac{41}{5} \right\rfloor+\left\lfloor\frac{41}{10} \right\rfloor
        +\left\lfloor\frac{41}{15} \right\rfloor
        -\left\lfloor\frac{41}{30} \right\rfloor
        +\left[\frac{41}{5}\frac{\widetilde{1}}{3}\right]
        -\left[
        \frac{41}{10}\frac{\widetilde{1}}{3}\right]\\
       & {\ \ }-\left[\frac{\widetilde{41}}{5} \right]
        +\left[\frac{41}{2}\frac{\widetilde{1}}{5} \right]
        +\left[\frac{41}{3}\frac{\widetilde{1}}{5} \right]
        -\left[\frac{41}{6}\frac{\widetilde{1}}{5} \right]
        +\left[\frac{\widetilde{41}}{15} \right]
        -\left[\frac{41}{2}\frac{\widetilde{1}}{15} \right]\\
       & =41-20
        -13+6
        -\left[\frac{41+2}{3} \right]
        +\left\lfloor\frac{
           \left\lfloor\frac{41}{2}\right\rfloor
            +3-2
         }{3}\right\rfloor\\
       & {\ \ }-8+4
        +2
        -1
        +\left\lfloor\frac{
           \left\lfloor\frac{41}{5}\right\rfloor
            +3-2
         }{3}\right\rfloor
        -\left\lfloor\frac{
           \left\lfloor\frac{41}{10}\right\rfloor
            +3-1
         }{3}\right\rfloor\\
       & {\ \ }-\left[\frac{41+2}{5} \right]
        +\left\lfloor\frac{
           \left\lfloor\frac{41}{2}\right\rfloor
            +5-3
         }{5}\right\rfloor
        +\left\lfloor\frac{
           \left\lfloor\frac{41}{3}\right\rfloor
            +5-1
         }{5}\right\rfloor
        -\left\lfloor\frac{
           \left\lfloor\frac{41}{6}\right\rfloor
            +5-3
         }{5}\right\rfloor
        +\left[\frac{41+2}{15} \right]
        -\left\lfloor\frac{
           \left\lfloor\frac{41}{2}\right\rfloor
            +15-14
         }{15}\right\rfloor\\
       & =14
        -14
        +\left\lfloor\frac{21}{3}\right\rfloor\\
       & {\ \ }-3
        +\left\lfloor\frac{9}{3}\right\rfloor
        -\left\lfloor\frac{6}{3}\right\rfloor\\
       & {\ \ }-8
        +\left\lfloor\frac{22}{5}\right\rfloor
        +\left\lfloor\frac{17}{5}\right\rfloor
        -\left\lfloor\frac{8}{5}\right\rfloor
        +\left[\frac{41+2}{15} \right]
        -\left\lfloor\frac{21}{15}\right\rfloor\\
       & =7+(-3+3-2)+(-8+4+3-1+2-1)\\
       & =4. \\
         \end{array}
\end{displaymath}
It is the same as before. From $D(\sqrt{41})=2$ for $p_i\leq 5$:
(3,5) and (5,7), $D_1(n)=0(n+2=43<p_4^2=49)$. So
$D(41)=D_0(41)+D(\sqrt{41})-D_1=4+2-0=6$.

The set: $Z=\{3,5,11,17,29,41\}, Z'=Z+2=\{5,7,13,19,31,43\}$. The
six twin primes are
$(3,5),(5,7),(11,13),(17,19),(29,31),(41,43)$.\qed
\end{example}

\begin{definition}
\label{def:2} The items of $Z_k\mod p_i= 0$, or the multiples of
$p_i$, have,
\begin{equation}
  \label{equ:d2}
   \begin{array}{lr}
    S(m,p_i\parallel0):=\sum\limits_{Z_k\mod  p_i=0}1
    =m\left[\frac{1}{p_i}\right]
    =\left\lfloor\frac{m}{p_i}\right\rfloor.
 \end{array}
\end{equation}
\end{definition}

So the items of  $Z_k\mod p_j\not= 0, p_j-2$ have,

\begin{equation}
  \label{equ:d22}
   \begin{array}{rclll}
    S(m,p_j\nparallel 0,p_j-2) & := &
    \sum\limits_{Z_k\mod p_i\ne0,p_j-2}1
     =
    m\left[1-\frac{1}{p_j}-\frac{\widetilde{1}}{p_j}\right].
 \end{array}
\end{equation}

\begin{definition}
\label{def:3} Let $X=\{X_1,X_2,\cdots,X_t\}$ be an (any) integer
set, then
\begin{equation}
  \label{equ:d3}
   \!\!\!\!\!\!\!\!
   \begin{array}{lr}
    t\left[1-\frac{1'}{p_i}-\frac{1''}{p_i}\right]
    :=\sum\limits_{    X_k\mod p_i\ne0,p_i-2 }1
    =t-\left[\frac{t+\phi'_i}{p_i}\right]-\left[\frac{t+\phi''_i}{p_i}\right]
      \geq0
 \end{array}
\end{equation}
is the number of items left after deleted the items of $X_k\mod
p_i=0,p_i-2$.
\end{definition}

Because
$m\left[1-\frac{1}{p_i}-\frac{\widetilde{1}}{p_i}\right]\left[1-\frac{1}{p_j}-\frac{\widetilde{1}}{p_j}\right]
\ne
\left(m\left[1-\frac{1}{p_i}-\frac{\widetilde{1}}{p_i}\right]\right)
\left[1-\frac{1}{p_j}-\frac{\widetilde{1}}{p_j}\right]$, We can
express it as
\begin{equation}
  \label{equ:d32}
 \begin{array}{rl}
    S(m,p_i\nparallel 0,p_i-2; p_j\nparallel 0,p_j-2)
    &=m\left[1-\frac{1}{p_i}-\frac{\widetilde{1}}{p_i}\right]\left[1-\frac{1}{p_j}-\frac{\widetilde{1}}{p_j}\right]\\
    &=\left(m\left[1-\frac{1}{p_i}-\frac{\widetilde{1}}{p_i}\right]\right)\left[1-\frac{1'}{p_j}-\frac{1''}{p_j}\right],
 \end{array}
\end{equation}
is the number of items left when we first delete those $Z_k\mod
p_i=0,p_i-2$ from set $Z$, and then delete those $X_{k'}\mod
p_j=0,p_j-2$ from set $X=\{Z, X_{k'}\mod p_i\ne 0,p_i-2,
k'=1,2,\cdots,
m\left[1-\frac{1}{p_i}-\frac{\widetilde{1}}{p_i}\right]\}$, where
set $X$ is no longer an arithmetic sequence. In general,

\begin{equation}
  \label{equ:d33}
 \begin{array}{lr}
    m\prod_{i=1}^{i_m}\limits\left[1-\frac{1}{p_i}-\frac{\widetilde{1}}{p_i}\right]
      \left[1-\frac{1}{p_j}-\frac{\widetilde{1}}{p_j}\right]
    =\left(m\prod_{i=1}^{i_m}\limits\left[1-\frac{1}{p_i}-\frac{\widetilde{1}}{p_i}\right]\right)
    \left[1-\frac{1'}{p_j}-\frac{1''}{p_j}\right].
 \end{array}
\end{equation}

\section{ Some property}
\label{sect:Property}

\begin{equation}
  \label{equ:p1}\left\{
   \begin{array}{lr}
    \left[1-\frac{1}{p_i}\right]\left[1-\frac{1}{p_j}\right]
    =\left[1-\frac{1}{p_j}\right]\left[1-\frac{1}{p_i}\right]\\
    \left[1-\frac{1}{p_i}-\frac{\widetilde{1}}{p_i}\right]\left[1-\frac{1}{p_j}-\frac{\widetilde{1}}{p_j}\right]
    =\left[1-\frac{1}{p_j}-\frac{\widetilde{1}}{p_j}\right]\left[1-\frac{1}{p_i}-\frac{\widetilde{1}}{p_i}\right].
 \end{array}\right.
\end{equation}

\begin{equation}
  \label{equ:p2}
   \begin{array}{rl}
      \left[\frac{m}{p_i}\right] & =\left[\frac{m_1+m_2}{p_i}\right]
         =\left[\frac{m_1}{p_i}\right]+\left[\frac{m_2}{p_i}\right]+\left[\frac{m_1
                \mod p_i+m_2 \mod p_i}{p_i}\right].
\end{array}
\end{equation}

\begin{equation}
  \label{equ:p3}
   \begin{array}{lr}
      m\left[1-\frac{1}{p_i}-\frac{\widetilde{1}}{p_i}\right]
      =ap_i\left(1-\frac{2}{p_i}\right)+b\left[1-\frac{1}{p_i}-\frac{\widetilde{1}}{p_i}\right]
      \mbox{\quad for $m=ap_i+b$}.
\end{array}
\end{equation}

\begin{equation}
  \label{equ:p4}
   \begin{array}{lr}
     -2\leq(m_1+m_2)\left[1-\frac{1}{p_i}-\frac{\widetilde{1}}{p_i}\right]
    -m_1\left[1-\frac{1}{p_i}-\frac{\widetilde{1}}{p_i}\right]
    -m_2\left[1-\frac{1}{p_i}-\frac{\widetilde{1}}{p_i}\right]\leq 1.\\
\end{array}
\end{equation}

\begin{proof} Let $\alpha=m_1\mod p_i,\beta=m_2\mod p_i,$
\begin{equation}
   \begin{array}{lr}
    \delta_{12}=(m_1+m_2)\left[1-\frac{1}{p_i}-\frac{\widetilde{1}}{p_i}\right]
    -m_1\left[1-\frac{1}{p_i}-\frac{\widetilde{1}}{p_i}\right]
    -m_2\left[1-\frac{1}{p_i}-\frac{\widetilde{1}}{p_i}\right]\nonumber\\
   =m_1+m_2-\left[\frac{m_1+m_2}{p_i}\right]-\left[\frac{m_1+m_2+2}{p_i}\right]
   -m_1+\left[\frac{m_1}{p_i}\right]+\left[\frac{m_1+2}{p_i}\right]\\
   \ \ -m_2+\left[\frac{m_2}{p_i}\right]+\left[\frac{m_2+2}{p_i}\right]\nonumber\\
   =-\left[\frac{\alpha+\beta}{p_i}\right]-\left[\frac{\alpha+\beta+2}{p_i}\right]
   +\left[\frac{\alpha+2}{p_i}\right]
   +\left[\frac{\beta+2}{p_i}\right]\nonumber.\\
 \end{array}
 \end{equation}

The minimum: $ min(\delta_{12})=-1-1=-2$, when $\alpha+\beta\geq
p_i$, $\alpha+2 <p_i$, $\beta+2 <p_i$, and $\alpha+\beta+2 \geq
p_i$.

The maximum: $ max(\delta_{12})=-0-1+1+1=1$, when $\alpha+\beta<
p_i$, $\alpha+2 \geq p_i$, $\beta+2 \geq p_i$.
 \end{proof}

We can represent it as
\begin{equation}
\begin{array}{lll}
  \label{equ:p5}
   &\ & (m_1+m_2)\left[1-\frac{1}{p_i}-\frac{\widetilde{1}}{p_i}\right]
     =  m_1+m_2-\left[\frac{m_1+m_2}{p_i}\right]-\left[\frac{m_1+m_2+2}{p_i}\right]\\
    & = & m_1-\left[\frac{m_1}{p_i}\right]-\left[\frac{m_1+2}{p_i}\right]
      +m_2-\left[\frac{m_1 \mod p_i+m_2}{p_i}\right]-\left[\frac{m_2+(m_1+2)\mod p_i}{p_i}\right]\\
    & = & m_1\left[1-\frac{1}{p_i}-\frac{\widetilde{1}}{p_i}\right]
    +m_2\left[1-\frac{1'}{p_i}-\frac{1''}{p_i}\right].\\
\end{array}
\end{equation}
$m_2\left[\frac{1'}{p_i}\right]$ will delete the items of $X_k
\mod p_i=0$ in the sequence of $X=\{m_1+1,m_1+2, \cdots,
m_1+m_2\}$, and $m_2\left[\frac{1''}{p_i}\right]$ will delete the
items of $X'_k \mod p_i=0$ in the sequence of $X'=\{X+2\}$ or the
items of $X_k \mod p_i=p_i-2=p_i-2$ in set $X$.

For any $m_2$, $0\leq
m_2\left[1-\frac{1'}{p_i}-\frac{1''}{p_i}\right] \leq m_2$, so,
\begin{equation}
  \label{equ:p6}
   \begin{array}{lr}
    %    0\leq(m+1)\left[1-\frac{1}{p_i}-\frac{\widetilde{1}}{p_i}\right]
    %    - m\left[1-\frac{1}{p_i}-\frac{\widetilde{1}}{p_i}\right]
    %    %=(1)\left[1-\frac{1'}{p_i}-\frac{1''}{p_i}\right]
    %    \leq    1\\
    m'\left[1-\frac{1}{p_i}-\frac{\widetilde{1}}{p_i}\right]
    \geq m\left[1-\frac{1}{p_i}-\frac{\widetilde{1}}{p_i}\right]
    \mbox{\quad for \ $m'\geq m$}.
 \end{array}
\end{equation}

%%%%%%%%%%%%%%%%%%%%%%%%%%%%%%%%%%%%%%%%%%%%%%%%%%%%%%%%%%%%%%%%%%%%%%%
\section{ Some lemma}
\label{sect:Lemma}

\begin{lemma} \label{lem:1}
For $m\geq p_i$,
\begin{equation}
  \label{equ:L1}
  \begin{array}{l}
     m\left[1-\frac{1}{p_i}-\frac{\widetilde{1}}{p_i}\right]
     \geq \left\lceil m\left(1-\frac{3}{p_i}\right)\right\rceil.\\
  \end{array}
\end{equation}
\end{lemma}

\begin{proof}
\begin{displaymath}
   \begin{array}{rl} left
  & =\left\lceil m-\left\lfloor\frac{m}{p_i}\right\rfloor-\left\lfloor\frac{m+\theta_i}{p_i}\right\rfloor\right\rceil
   \geq\left\lceil
   m-\left\lfloor\frac{m}{p_i}\right\rfloor-\left\lfloor\frac{m}{p_i}\right\rfloor-1\right\rceil\\
   & \geq\left\lceil m-\frac{2m}{p_i}-\frac{m}{p_i}\right\rceil
    = right.
\end{array}
\end{displaymath}
\end{proof}

\begin{lemma} \label{lem:2}

\begin{equation}
 \label{equ:L2}
 \left\{
 \begin{array}{lll}
    -1 \leq & m\left[\frac{1}{p_j}\right]\left[1-\frac{1}{p_i}-\frac{\widetilde{1}}{p_i}\right]
    -\left\lfloor\frac{m}{p_j}\right\rfloor\left[1-\frac{1}{p_i}-\frac{\widetilde{1}}{p_i}\right] & \leq1\\
    -1 \leq & m\left[\frac{\widetilde{1}}{p_j}\right]\left[1-\frac{1}{p_i}-\frac{\widetilde{1}}{p_i}\right]
    -\left\lfloor\frac{m+2}{p_j}\right\rfloor\left[1-\frac{1}{p_i}-\frac{\widetilde{1}}{p_i}\right] & \leq1.\\
\end{array}\right.
\end{equation}
\end{lemma}

\begin{proof}
\begin{equation}\nonumber
 \begin{array}{ll}
   \varepsilon_1 &=m\left[\frac{1}{p_j}\right]\left[1-\frac{1}{p_i}-\frac{\widetilde{1}}{p_i}\right]
    -\left\lfloor\frac{m}{p_j}\right\rfloor\left[1-\frac{1}{p_i}-\frac{\widetilde{1}}{p_i}\right]\\
    &=\left\lfloor\frac{m}{p_j}\right\rfloor-\left\lfloor\frac{m}{p_ip_j}\right\rfloor
      -\left\lfloor\frac{m+\theta_{j,i}}{p_ip_j}\right\rfloor
      -\left\lfloor\frac{m}{p_j}\right\rfloor+\left\lfloor\frac{m}{p_ip_j}\right\rfloor
      +\left\lfloor\frac{\left\lfloor\frac{m}{p_j}\right\rfloor+2}{p_i}\right\rfloor\\
    &= -\left\lfloor\frac{m+\theta_{j,i}}{p_ip_j}\right\rfloor
      +\left\lfloor\frac{m+2p_j}{p_ip_j}\right\rfloor\\
    &= -\left\lfloor\frac{m\mod p_ip_j+\theta_{j,i}}{p_ip_j}\right\rfloor
      +\left\lfloor\frac{m\mod p_ip_j+2p_j}{p_ip_j}\right\rfloor,\\
\end{array}
\end{equation}
so that $-1\leq \varepsilon_1\leq1$.
\begin{equation}\nonumber
 \begin{array}{ll}
   \varepsilon_2= m\left[\frac{\widetilde{1}}{p_j}\right]\left[1-\frac{1}{p_i}-\frac{\widetilde{1}}{p_i}\right]
    -\left\lfloor\frac{m+2}{p_j}\right\rfloor\left[1-\frac{1}{p_i}-\frac{\widetilde{1}}{p_i}\right]\\
    = \left\lfloor\frac{m+2}{p_j}\right\rfloor
        -\left\lfloor\frac{m+\theta_{i,j}}{p_ip_j}\right\rfloor
        -\left\lfloor\frac{m+2}{p_ip_j}\right\rfloor
    -\left\lfloor\frac{m+2}{p_j}\right\rfloor
        +\left\lfloor\left\lfloor\frac{m+2}{p_j}\right\rfloor\frac{1}{p_i}\right\rfloor
        +\left\lfloor\frac{\left\lfloor\frac{m+2}{p_j}\right\rfloor+2}{p_i}\right\rfloor\\
    =   -\left\lfloor\frac{m\mod p_ip_j+\theta_{i,j}}{p_ip_j}\right\rfloor
        +\left\lfloor\frac{m\mod p_ip_j+2+2p_j}{p_ip_j}\right\rfloor,\\
\end{array}
\end{equation}
so that $-1\leq \varepsilon_2\leq1$.
\end{proof}

\begin{lemma} \label{lem:3}

For $m\geq p_j^2$, $p_j>p_i\geq 3$,
\begin{equation}
 \label{equ:L3}
   \begin{array}{lr}
    m\left[1-\frac{1}{p_i}-\frac{\widetilde{1}}{p_i}\right]
    \left[1-\frac{1}{p_j}-\frac{\widetilde{1}}{p_j}\right]
    \geq \left\lceil m\left(1-\frac{3}{p_j}\right)\right\rceil
    \left[1-\frac{1}{p_i}-\frac{\widetilde{1}}{p_i}\right].\\
\end{array}
\end{equation}

\begin{equation}
 \label{equ:L32}
   \begin{array}{lr}
    m\left[1-\frac{1}{p_i}\right]
    \left[1-\frac{1}{p_j}-\frac{\widetilde{1}}{p_j}\right]
    \geq \left\lceil m\left(1-\frac{3}{p_j}\right)\right\rceil
    \left[1-\frac{1}{p_i}\right]\mbox{\quad for $p_i=2$}.
\end{array}
\end{equation}

\end{lemma}

%%%%%%%%%%%%%%%%%%%
\begin{proof}[Proof of (\ref{equ:L3}).]
Let $m=sp_ip_j+t$, $t=ap_j+b$, $0\leq a\leq (p_i-1),\ 0\leq b\leq
(p_j-1)$, Because $m\geq p_j^2,\ p_j>p_i$, so  $s\geq 1$. From Eq.
(\ref{equ:p2}), (\ref{equ:p3}),

\begin{displaymath}
   \begin{array}{rl}
    \varepsilon & =m\left[1-\frac{1}{p_i}-\frac{\widetilde{1}}{p_i}\right]
         \left[1-\frac{1}{p_j}-\frac{\widetilde{1}}{p_j}\right]
         -\left\lceil m\left(1-\frac{3}{p_j}\right)\right\rceil
         \left[1-\frac{1}{p_i}-\frac{\widetilde{1}}{p_i}\right]\\
     & =sp_ip_j\left(1-\frac{2}{p_i}\right)\left(1-\frac{2}{p_j}\right)
        -sp_ip_j\left(1-\frac{2}{p_i}\right)\left(1-\frac{3}{p_j}\right)\\
        &\ \ +t-\left\lfloor\frac{t}{p_i}\right\rfloor-\left\lfloor\frac{t+2}{p_i}\right\rfloor
         -\left\lfloor\frac{t}{p_j}\right\rfloor+\left\lfloor\frac{t}{p_ip_j}\right\rfloor
         +\left\lfloor\frac{t+\theta_{j,i}}{p_ip_j}\right\rfloor
         -\left\lfloor\frac{t+2}{p_j}\right\rfloor+\left\lfloor\frac{t+\theta_{i,j}}{p_ip_j}\right\rfloor
         +\left\lfloor\frac{t+2}{p_ip_j}\right\rfloor\\
         &\ \ -\left(t-\left\lfloor\frac{3t}{p_j}\right\rfloor\right)
         +\left\lfloor\frac{t-\left\lfloor\frac{3t}{p_j}\right\rfloor}{p_i}\right\rfloor
         +\left\lfloor\frac{t+2-\left\lfloor\frac{3t}{p_j}\right\rfloor}{p_i}\right\rfloor\\
     & =s(p_i-2)\\
        &\ \ +\left\lfloor\frac{3t}{p_j}\right\rfloor
         -\left\lfloor\frac{t}{p_j}\right\rfloor-\left\lfloor\frac{t+2}{p_j}\right\rfloor
         +\left\lfloor\frac{t}{p_ip_j}\right\rfloor
         +\left\lfloor\frac{t+\theta_{j,i}}{p_ip_j}\right\rfloor
         +\left\lfloor\frac{t+\theta_{i,j}}{p_ip_j}\right\rfloor
         +\left\lfloor\frac{t+2}{p_ip_j}\right\rfloor\\
         &\ \ -\left\lfloor\frac{\left\lfloor\frac{3t}{p_j}\right\rfloor}{p_i}\right\rfloor
         -\varepsilon_1
          -\left\lfloor\frac{\left\lfloor\frac{3t}{p_j}\right\rfloor}{p_i}\right\rfloor
         -\varepsilon_2\\
     & =s(p_i-2)\\
        &\ \ +3a+\left\lfloor\frac{3b}{p_j}\right\rfloor
         -a-a-\left\lfloor\frac{b+2}{p_j}\right\rfloor
         +0
         +\left\lfloor\frac{ap_j+b+\theta_{j,i}}{p_ip_j}\right\rfloor
         +\left\lfloor\frac{ap_j+b+\theta_{i,j}}{p_ip_j}\right\rfloor
         +\left\lfloor\frac{ap_j+b+2}{p_ip_j}\right\rfloor\\
         &\ \ -2\left\lfloor\frac{3a+\left\lfloor\frac{3b}{p_j}\right\rfloor}{p_i}\right\rfloor
         -\varepsilon_1     -\varepsilon_2\\
     & =s(p_i-2)+a+\left\lfloor\frac{3b}{p_j}\right\rfloor
         +\left\lfloor\frac{ap_j+b+\theta_{j,i}}{p_ip_j}\right\rfloor
         +\left\lfloor\frac{ap_j+b+\theta_{i,j}}{p_ip_j}\right\rfloor
         +\left\lfloor\frac{ap_j+b+2}{p_ip_j}\right\rfloor\\
         &\ \
         -\varepsilon_1         -\varepsilon_2-2\varepsilon_3-\varepsilon_4,\\
\end{array}
\end{displaymath}
where
\begin{displaymath}
   \left\{
   \begin{array}{rl}
     \varepsilon_1 & =
         \left\lfloor\frac{\left\lfloor\frac{3t}{p_j}\right\rfloor\mod p_i
         +\left(t-\left\lfloor\frac{3t}{p_j}\right\rfloor\right)\mod p_i}
         {p_i}\right\rfloor\leq 1\\
     \varepsilon_2 & =
         \left\lfloor\frac{\left\lfloor\frac{3t}{p_j}\right\rfloor\mod p_i
         +\left(t+2-\left\lfloor\frac{3t}{p_j}\right\rfloor\right)\mod p_i}
         {p_i}\right\rfloor\leq 1\\
     \varepsilon_3 & =
         \left\lfloor\frac{3a+\left\lfloor\frac{3b}{p_j}\right\rfloor}{p_i}\right\rfloor\leq 2\\
     \varepsilon_4 & =
         \left\lfloor\frac{b+2}{p_j}\right\rfloor\leq 1,\\
\end{array}\right.
\end{displaymath}

\begin{displaymath}
     \Delta
     \varepsilon=\varepsilon_1+\varepsilon_2+2\varepsilon_3+\varepsilon_4\leq7.
\end{displaymath}

  \begin{enumerate}
  \item  If $p_i\geq 11$, then $\varepsilon\geq
  s(11-2)-\Delta \varepsilon\geq 9-7>0.$

  \item  If $p_i=7$, then
    \begin{enumerate}
    \item  If $a\geq2$, then $\varepsilon\geq s(p_i-2)+a-\Delta
    \varepsilon\geq 5+2-7=0.$

    \item  If $a\leq1$, then $\varepsilon_3  =
         \left\lfloor\frac{3a+\left\lfloor\frac{3b}{p_j}\right\rfloor}{p_i}\right\rfloor
         \leq \left\lfloor\frac{3+2}{7}\right\rfloor=0
    $, $\Delta \varepsilon\leq3$, so $\varepsilon\geq s(p_i-2)-\Delta
    \varepsilon\geq 5-3>0.$
    \end{enumerate}

  \item  If $p_i=5$, then
    \begin{enumerate}

    \item  If $a\geq4$, then $\varepsilon\geq s(p_i-2)+a-\Delta
    \varepsilon\geq 3+4-7=0.$

    \item  If $a=0$, then $\varepsilon_3  =0 $, $\Delta
    \varepsilon\leq3$, so $\varepsilon\geq s(p_i-2)-\Delta
    \varepsilon\geq 3-3=0.$

    \item  If $a=1$: if $\left\lfloor\frac{3b}{p_j}\right]\geq1$ then
    $\varepsilon_3 \leq1 $, $\Delta \varepsilon\leq5$, so
    $\varepsilon\geq
    s(p_i-2)+a+\left\lfloor\frac{3b}{p_j}\right]-\Delta
    \varepsilon\geq 3+1+1-5=0$; else for
    $\left\lfloor\frac{3b}{p_j}\right]=0$ then $\varepsilon_3 =0 $,
    $\Delta \varepsilon\leq3$, so $\varepsilon\geq s(p_i-2)+a-\Delta
    \varepsilon\geq 3+1-3>0$.

    \item  If $a=2$, then $\varepsilon_3  =1 $, $\Delta
    \varepsilon\leq5$, so $\varepsilon\geq s(p_i-2)+a-\Delta
    \varepsilon\geq 3+2-5=0.$

    \item  If $a=3$: if $\left\lfloor\frac{3b}{p_j}\right]\geq1$ then
    $\varepsilon\geq
    s(p_i-2)+a+\left\lfloor\frac{3b}{p_j}\right]-\Delta
    \varepsilon\geq 3+3+1-7=0$; else for
    $\left\lfloor\frac{3b}{p_j}\right]=0$ then $\varepsilon_3 =
    \left\lfloor\frac{3a+\left\lfloor\frac{3b}{p_j}\right\rfloor}{p_i}\right\rfloor=1
    $, $\Delta \varepsilon\leq5$, so $\varepsilon\geq
    s(p_i-2)+a-\Delta \varepsilon\geq 3+3-5>0$.
    \end{enumerate}

  \item  If $p_i=3$, then

  \begin{displaymath}
   \left\{
   \begin{array}{rl}
     a-2\varepsilon_3 & =a -2\left\lfloor\frac{3a+\left\lfloor\frac{3b}{p_j}\right\rfloor}{3}\right\rfloor
     =a-2a=-a\\
     \varepsilon_1 & =
         \left\lfloor\frac{\left(3a+\left\lfloor\frac{3b}{p_j}\right\rfloor\right)\mod 3
         +\left(ap_j+b-3a-\left\lfloor\frac{3b}{p_j}\right\rfloor\right)\mod 3}
         {3}\right\rfloor\\
   & =   \left\lfloor\frac{\left(\left\lfloor\frac{3b}{p_j}\right\rfloor\right)\mod 3
         +\left(ap_j+b-\left\lfloor\frac{3b}{p_j}\right\rfloor\right)\mod 3}
         {3}\right\rfloor\\
     \varepsilon_2 & =
         \left\lfloor\frac{\left(3a+\left\lfloor\frac{3b}{p_j}\right\rfloor\right)\mod 3
         +\left(ap_j+b+2-3a-\left\lfloor\frac{3b}{p_j}\right\rfloor\right)\mod 3}
         {3}\right\rfloor\\
   & =   \left\lfloor\frac{\left(\left\lfloor\frac{3b}{p_j}\right\rfloor\right)\mod 3
         +\left(ap_j+b+2-\left\lfloor\frac{3b}{p_j}\right\rfloor\right)\mod 3}
         {3}\right\rfloor,\\
  \end{array}\right.
  \end{displaymath}

  \begin{displaymath}
  \label{equ:e0}
   \begin{array}{rl}
    \varepsilon & =s-a+\left\lfloor\frac{3b}{p_j}\right\rfloor
         +\left\lfloor\frac{ap_j+b+\theta_{j,i}}{p_ip_j}\right\rfloor
         +\left\lfloor\frac{ap_j+b+\theta_{i,j}}{p_ip_j}\right\rfloor
         +\left\lfloor\frac{ap_j+b+2}{p_ip_j}\right\rfloor\\
         &\ \
         -\varepsilon_1         -\varepsilon_2-\varepsilon_4.\\
  \end{array}
  \end{displaymath}

    \begin{enumerate}
    \item  If $\left\lfloor\frac{3b}{p_j}\right\rfloor=0 $, then
    $\varepsilon_1=\varepsilon_2=0$, $b<p_j/3,\
    b+2<p_j/3+2=(p_j+6)/3<p_j$, so $\varepsilon_4=0$. If $a\leq1$
    then $\varepsilon\geq s-a\geq0$. Else for $a=2 (<p_i)$, because
    $\theta_{j,i}=p_ip_j-\lambda p_j\geq p_j(\ 1\leq \lambda\leq
    p_i-1)$,
    $\left\lfloor\frac{ap_j+b+\theta_{j,i}}{p_ip_j}\right\rfloor\geq
    \left\lfloor\frac{2p_j+p_j}{3p_j}\right\rfloor=1$. So
    $\varepsilon\geq s-a+1\geq0$.

    \item If $\left\lfloor\frac{3b}{p_j}\right\rfloor\geq1 $, then
    $
      \label{equ:sp35}
       \begin{array}{rl}
        \varepsilon' =s+\left\lfloor\frac{ap_j+b+\theta_{i,j}}{3p_j}\right\rfloor
        -\varepsilon_4\geq1.\\
    \end{array}$

    \begin{proof}
      \begin{itemize}

      \item  If $s\geq 2$ then $\varepsilon'\geq 1$.

      \item  If $s=1$, i.e.,
      $s=\left\lfloor\frac{m}{p_ip_j}\right\rfloor\geq\left\lfloor\frac{p_j^2}{3p_j}\right\rfloor=1$,
      so $p_j=5$. From $5^2=3\cdot5+2\cdot5$, we have $a=2$. Let us
      consider $\theta_{i,j}=p_ip_j-\lambda_{i,j}$,
      $\lambda_{i,j}=\lambda p_i=3\lambda$, with the condition
      $(3\lambda+2)\mod 5=0$ i.e. $ \lambda=1$,
      $\theta_{i,j}=p_ip_j-3\lambda=15-3=12$.
      $\left\lfloor\frac{ap_j+b+\theta_{i,j}}{p_ip_j}\right\rfloor \geq
      \left\lfloor\frac{2\cdot 5+12}{3\cdot5}\right\rfloor=1$,
      Therefore,
      \begin{displaymath}
        \label{equ:e1}
         \begin{array}{rl}
          \varepsilon' =s+\left\lfloor\frac{ap_j+b+\theta_{i,j}}{3p_j}\right\rfloor
          -\varepsilon_4\geq 1 \quad\mbox{if}\quad s=1.\\
      \end{array}
      \end{displaymath}
      \end{itemize}
    \end{proof}

      Besides, $
         \begin{array}{rl}
          \varepsilon'' & =\left\lfloor\frac{3b}{p_j}\right\rfloor
               -\varepsilon_1 -\varepsilon_2\geq0.\\
      \end{array}$

      \begin{proof}
      \begin{itemize}
      \item  if $\left\lfloor\frac{3b}{p_j}\right\rfloor=2 $, then
      $\varepsilon''\geq 0$.

      \item   if $\left\lfloor\frac{3b}{p_j}\right\rfloor=1 $, because,

      \begin{displaymath}
         \begin{array}{l}
          \varepsilon_1+\varepsilon_2=\left\lfloor\frac{1
               +\left(ap_j+b-1\right)\mod 3} {3}\right\rfloor
               +\left\lfloor\frac{1
               +\left(ap_j+b+1\right)\mod 3} {3}\right\rfloor\\
               =\left\{
         \begin{array}{rl}
           1+0=1\quad\mbox{if}\quad (ap_j+b)\mod 3=0\\
           0+1=1\quad\mbox{if}\quad (ap_j+b)\mod 3=1\\
           0+0=0\quad\mbox{if}\quad (ap_j+b)\mod 3=2.\\
      \end{array}\right.
      \end{array}
      \end{displaymath}

      So $\varepsilon_1+\varepsilon_2\leq1$ and
      $
      \label{equ:e2}
         \begin{array}{rl}
          \varepsilon'' & =\left\lfloor\frac{3b}{p_j}\right\rfloor
         -\varepsilon_1 -\varepsilon_2\geq0\\
      \end{array}.$
      \end{itemize}
      \end{proof}

    Therefore, %From Eq. (\ref{equ:e0}), (\ref{equ:e1}), (\ref{equ:e2}),
    \begin{displaymath}
       \begin{array}{rl}
        \varepsilon & =\varepsilon'-a+\varepsilon''
             +\left\lfloor\frac{ap_j+b+\theta_{j,i}}{p_ip_j}\right\rfloor
             +\left\lfloor\frac{ap_j+b+2}{p_ip_j}\right\rfloor\\
       & \geq
       1-a+\left\lfloor\frac{ap_j+b+\theta_{j,i}}{p_ip_j}\right\rfloor.\\
    \end{array}
    \end{displaymath}

    If $a\leq 1$ then $\varepsilon\geq 1-a\geq0$.

    Else for $a=2 (<p_i)$, because $\theta_{j,i}=p_ip_j-\lambda p_j\geq
    p_j\ (1\leq\lambda\leq p_i-1)$,
    $\left\lfloor\frac{ap_j+b+\theta_{j,i}}{p_ip_j}\right\rfloor\geq
    \left\lfloor\frac{2p_j+p_j}{3p_j}\right\rfloor=1$. So
    $\varepsilon\geq 1-a+1\geq0$.

{\ }
    \end{enumerate}
  \end{enumerate}
In summary, $\varepsilon \geq 0$ for all $p_i<p_j\leq \sqrt{m}$.
\end{proof}

%%%%%%%%%%%%%%%%%%%%%%%%%%

\begin{proof}[Proof of (\ref{equ:L32}).]
when $p_i=2$, $\left[\frac{\widetilde{1}}{p_i}\right]\equiv0$,
$m=2sp_j+t,\ t=ap_j+b$, $0\leq a\leq 1$, $0\leq b\leq p_j-1$,
$s\geq 1$.

\begin{displaymath}
   \begin{array}{rl}
    \varepsilon & =m\left[1-\frac{1}{2}\right]
         \left[1-\frac{1}{p_j}-\frac{\widetilde{1}}{p_j}\right]
         -\left\lceil m\left(1-\frac{3}{p_j}\right)\right\rceil
         \left[1-\frac{1}{2}\right]\\
     & =2sp_j\left(1-\frac{1}{2}\right)\left(1-\frac{2}{p_j}\right)
        -2sp_j\left(1-\frac{1}{2}\right)\left(1-\frac{3}{p_j}\right)\\
        &\ \ +t-\left\lfloor\frac{t}{2}\right\rfloor
         -\left\lfloor\frac{t}{p_j}\right\rfloor+\left\lfloor\frac{t}{2p_j}\right\rfloor
         -\left\lfloor\frac{t+2}{p_j}\right\rfloor+\left\lfloor\frac{t+\theta_{i,j}}{2p_j}\right\rfloor
         \\
         &\ \ -\left(t-\left\lfloor\frac{3t}{p_j}\right\rfloor\right)
         +\left\lfloor\frac{t-\left\lfloor\frac{3t}{p_j}\right\rfloor}{2}\right\rfloor\\
     & =s+\left\lfloor\frac{3t}{p_j}\right\rfloor
         -\left\lfloor\frac{t}{p_j}\right\rfloor-\left\lfloor\frac{t+2}{p_j}\right\rfloor
         +0+\left\lfloor\frac{t+\theta_{i,j}}{2p_j}\right\rfloor
         -\left\lfloor\frac{\left\lfloor\frac{3t}{p_j}\right\rfloor}{2}\right\rfloor
         -\varepsilon_1\\
     & =s+3a+\left\lfloor\frac{3b}{p_j}\right\rfloor
         -a-a-\left\lfloor\frac{b+2}{p_j}\right\rfloor
         +\left\lfloor\frac{ap_j+b+\theta_{i,j}}{2p_j}\right\rfloor
         -\left\lfloor\frac{3a+\left\lfloor\frac{3b}{p_j}\right\rfloor}{2}\right\rfloor
         -\varepsilon_1\\
     & =s+\left\lfloor\frac{3b}{p_j}\right\rfloor
         +\left\lfloor\frac{ap_j+b+\theta_{i,j}}{2p_j}\right\rfloor
         -\varepsilon_1-\varepsilon_3-\varepsilon_4,\\
\end{array}
\end{displaymath}

where
\begin{displaymath}
   \left\{
   \begin{array}{rl}
     \varepsilon_1 & =
         \left\lfloor\frac{\left\lfloor3a+\frac{3b}{p_j}\right\rfloor\mod 2
         +\left(ap_j+b-3a-\left\lfloor\frac{3b}{p_j}\right\rfloor\right)\mod 2}
         {2}\right\rfloor\\
         & =
         \left\lfloor\frac{\left\lfloor a+\frac{3b}{p_j}\right\rfloor\mod 2
         +\left(ap_j+b-a-\left\lfloor\frac{3b}{p_j}\right\rfloor\right)\mod 2}
         {2}\right\rfloor\leq1\\
     \varepsilon_3 & =
         \left\lfloor\frac{a+\left\lfloor\frac{3b}{p_j}\right\rfloor}{2}\right\rfloor\leq1\\
     \varepsilon_4 & =
         \left\lfloor\frac{b+2}{p_j}\right\rfloor\leq1,\\
\end{array}\right.
\end{displaymath}

\begin{displaymath}
     \Delta\varepsilon
     =\varepsilon_1+\varepsilon_3+\varepsilon_4\leq 3.
\end{displaymath}

  \begin{enumerate}

  \item  If $\left\lfloor\frac{3b}{p_j}\right\rfloor=0$, then
    $\varepsilon_3=0$, $      \varepsilon_1  =
         \left\lfloor\frac{a\mod 2
         +\left(ap_j+b-a\right)\mod 2}
         {2}\right\rfloor:$

    If $a=0$, then $\varepsilon_1=0$, $\Delta\varepsilon
     =\varepsilon_4\leq 1$, $\varepsilon\geq
     s-\Delta\varepsilon\geq0$;

    Else for $a=1$, then $\varepsilon_1 =\left\lfloor\frac{1+b\mod 2}
    {2}\right\rfloor$. $\lambda_{i,j}=2p_j-2$,
    $\theta_{i,j}=2p_j-\lambda_{i,j}=2$.

    $\left\lfloor\frac{ap_j+b+\theta_{i,j}}{2p_j}\right\rfloor-\varepsilon_4
    =\left\lfloor\frac{p_j+b+2}{2p_j}\right\rfloor-\left\lfloor\frac{b+2}{p_j}\right\rfloor
    =\left\{
      \begin{array}{l}
          0-0=0\quad\mbox{if}\quad b+2<p_j\\
          1-1=0\quad\mbox{if}\quad b+2\geq p_j,\\
      \end{array}
    \right.$

    $\left\lfloor\frac{3b}{p_j}\right\rfloor-\varepsilon_3
    =\left\lfloor\frac{3b}{p_j}\right\rfloor
    -\left\lfloor\frac{a+\left\lfloor\frac{3b}{p_j}\right\rfloor}{2}\right\rfloor
    \geq\left\{
      \begin{array}{l}
          0-0=0\quad\mbox{if}\quad \left\lfloor\frac{3b}{p_j}\right\rfloor=0\\
          1-1=0\quad\mbox{if}\quad \left\lfloor\frac{3b}{p_j}\right\rfloor=1\\
          2-1=1\quad\mbox{if}\quad \left\lfloor\frac{3b}{p_j}\right\rfloor=2,\\
      \end{array}
    \right.$

    So $ \varepsilon =s+\left(\left\lfloor\frac{3b}{p_j}\right\rfloor
    -\varepsilon_3\right)
    +\left(\left\lfloor\frac{ap_j+b+\theta_{i,j}}{2p_j}\right\rfloor-\varepsilon_4\right)
         -\varepsilon_1
         \geq 1+0+0-1=0.$
\item
  If $\left\lfloor\frac{3b}{p_j}\right\rfloor=1$, then
  $\varepsilon\geq s+1+0-\varepsilon_1-\varepsilon_3-1 \geq
  1-\varepsilon_1-\varepsilon_3$:

    If $a=0$, then $\varepsilon_3=0$, $\varepsilon\geq
    1-\varepsilon_1-0\geq0$;

    Else for $a=1$, then $\varepsilon_1=0$, $\varepsilon\geq
    1-0-\varepsilon_3\geq0$.
\item
  If $\left\lfloor\frac{3b}{p_j}\right\rfloor=2$, then
$\varepsilon\geq s+2+0-\varepsilon_1-\varepsilon_3-\varepsilon_4
\geq 0$.

%  \end{enumerate}
   \end{enumerate}

In summary, $\varepsilon \geq 0$ for all $2=p_i<p_j\leq \sqrt{m}$.
\end{proof}

%%%%%%%%%%%%%%%%%%%%%%%%%%%%%%%%%%%

Eq. (\ref{equ:L3}) means that,
%for any $0\leq r_{p_i}\leq p_i-1$,
after deleted the items of $Z_k\mod p_j=0,p_j-2$ from $Z$, the
items of $Z_k\mod p_i\ne 0, p_i-2$ will have at least,
\begin{equation}
 \label{equ:L33}
   \begin{array}{l}
    S(m,p_i\nparallel 0,p_i-2; p_j\nparallel 0,p_j-2)\\
    = S(\lceil m(1-3/p_j)\rceil,p_i\nparallel 0,p_i-2)+\Delta(p_i\nparallel 0,p_i-2),
 \end{array}
\end{equation}
where $\Delta(p_i\nparallel 0,p_i-2)\geq0$ is the extra items of
$Z_k\mod p_i\ne0, p_i-2$ and $Z_k\mod p_j\ne 0,p_j-2$. Thus the
effect of $\left[1-\frac{1}{p_j}-\frac{\widetilde{1}}{p_j}\right]$
is that it constructs a new effective nature sequence with at
least $\lceil m(1-3p_j)\rceil$ items which satisfy the condition
$Z_k\mod p_j\ne0,p_j-2$.

This lemma means that, from equation (\ref{equ:L1}) and
(\ref{equ:p5}), we can let $
m\left[1-\frac{1}{p_j}-\frac{\widetilde{1}}{p_j}\right]=
     \left\lceil m\left(1-\frac{3}{p_j}\right)\right\rceil+t,
    t\geq 0$, and operated by
    $\left[1-\frac{1}{p_i}-\frac{\widetilde{1}}{p_i}\right]$ to
    get the items which have no multiples of $p_i$ and $p_j$ in $Z$ and$Z'$.

\begin{equation}
 \label{equ:L34}
   \begin{array}{lrclll}
    m\left[1-\frac{1}{p_i}-\frac{\widetilde{1}}{p_i}\right]
    \left[1-\frac{1}{p_j}-\frac{\widetilde{1}}{p_j}\right]
    =S(m,p_i\nparallel 0,p_i-2; p_j\nparallel 0,p_j-2)\\
    =\left( \left\lceil
    m\left(1-\frac{3}{p_j}\right)\right\rceil+t\right)
    \left[1-\frac{1}{p_i}-\frac{\widetilde{1}}{p_i}\right]\\
    =S(\lceil
    m(1-3/p_j)\rceil,p_i\nparallel 0,p_i-2)
    +\Delta (p_i\nparallel 0,p_i-2) \\
    =\left\lceil
    m\left(1-\frac{3}{p_j}\right)\right\rceil
    \left[1-\frac{1}{p_i}-\frac{\widetilde{1}}{p_i}\right]
    +t\left[1-\frac{1'}{p_i}-\frac{1''}{p_i}\right]\\
    \geq \left\lceil
    m\left(1-\frac{3}{p_j}\right)\right\rceil
    \left[1-\frac{1}{p_i}-\frac{\widetilde{1}}{p_i}\right].
 \end{array}
\end{equation}

%%%%%%%%%%%%%%%%%%%%%%%%%%%%%%%%%%%%%%%%%%%%%%%%%%%%%%%%%%%%%%%%%%%%%%%%%%%%%%%%%%%%%%%%%%%%%%%%%%%%

\begin{lemma} \label{lem:4}

\begin{equation}
  \label{equ:L4}
   \begin{array}{lr}
   m\left[1-\frac{1}{2}\right] \left[1-\frac{1}{3}-\frac{\widetilde{1}}{3}\right]
    \geq\left\lceil\frac{m}{6}\right\rceil-1.\\
 \end{array}
\end{equation}
\end{lemma}
\begin{proof}
For $p_i=2, p_j=3$, $\frac{\widetilde{1}}{2}=0$, let $m=6s+t$,
\begin{displaymath}
   \begin{array}{lr}
    m\left[1-\frac{1}{2}\right]
    \left[1-\frac{1}{3}-\frac{\widetilde{1}}{3}\right]
    = 6s\left[1-\frac{1}{2}\right]
    \left[1-\frac{1}{3}-\frac{\widetilde{1}}{3}\right]
    +t\left[1-\frac{1}{2}\right]
    \left[1-\frac{1}{3}-\frac{\widetilde{1}}{3}\right]\\
    \geq 6s\left(1-\frac{1}{2}\right)\left(1-\frac{2}{3}\right)
    =s=\left\lfloor\frac{m}{6}\right\rfloor
    \geq \left\lceil\frac{m}{6}\right\rceil-1.
 \end{array}
\end{displaymath}

For the residual class of modulo 6, $X=\{6s+1,6s+2,\cdots,6s+6\}$,
$X_k \mod 6=\{1, 2, 3, 4, 5, 6\}$, there are 4 elements $(2, 3, 4,
6)$ of multiples of 2 or 3. For the other elements $X_k \mod
6=\{1, 5\}$,
% which are coprime with 2 and 3,
there is at least one with $X_k\mod 6=5$, i.e., $X_k\mod 2\not
=0,X_k\mod 3\not =0$, and $X_k\mod 3\not = (3-2)=1$. The twin
primes have at least
$\left\lfloor\frac{m}{6}\right\rfloor\geq\left\lceil\frac{m}{6}\right\rceil-1$,
or
$m\left[1-\frac{1}{2}\right]\left[1-\frac{1}{3}-\frac{\widetilde{1}}{3}\right]
\geq\left\lceil\frac{m}{6}\right\rceil-1$.
\end{proof}

\begin{lemma} \label{lem:5}

For $ i=1, 2, \cdots, i_m, j$,
\begin{equation}
 \label{equ:L5}
   \begin{array}{lr}
    m\prod_{i=1}^{i_m}\limits\left[1-\frac{1}{p_i}-\frac{\widetilde{1}}{p_i}\right]
    \left[1-\frac{1}{p_j}-\frac{\widetilde{1}}{p_j}\right]
    \geq \left\lceil m\left(1-\frac{3}{p_j}\right)\right\rceil
    \prod_{i=1}^{i_m}\left[1-\frac{1}{p_i}-\frac{\widetilde{1}}{p_i}\right],\\
\end{array}
\end{equation}
where $\frac{\widetilde{1}}{p_i}\equiv 0$ for $p_i=2$.
\end{lemma}

\begin{proof}
From equation
(\ref{equ:d22}),(\ref{equ:L1}),(\ref{equ:L3})-(\ref{equ:L34}), For
$i_1\ne i_2$, we have $S(m,p_{i_1}\;\nparallel
0,p_{i_1}-2;p_j\;\nparallel 0,p_j-2) \geq S(\lceil
m(1-3/p_j)\rceil,p_{i_1}\;\nparallel 0,p_{i_1}-2)$,
$S(m,p_{i_2}\;\nparallel 0,p_{i_2}-2;p_j\;\nparallel 0,p_j-2) \geq
S(\lceil m(1-3/p_j)\rceil,p_{i_2}\;\nparallel 0,p_{i_2}-2)$, so
that $S(m,p_{i_1}\;\nparallel 0,p_{i_1}-2; p_{i_2}\;\nparallel
0,p_{i_2}-2; p_j\;\nparallel 0,p_j-2) \geq S(\lceil
m(1-3/p_j)\rceil,p_{i_1}\;\nparallel 0,p_{i_1}-2;
p_{i_2}\;\nparallel 0,p_{i_2}-2)$. $\cdots\cdots$.
$S(m,p_{i_1}\;\nparallel 0,p_{i_1}-2; p_{i_2}\;\nparallel
0,p_{i_2}-2;\cdots,p_{i_m}\;\nparallel 0,p_{i_m}-2;
p_j\;\nparallel 0,p_j-2) \geq S(\lceil
m(1-3/p_j)\rceil,p_{i_1}\;\nparallel 0,p_{i_1}-2;
p_{i_2}\;\nparallel 0,p_{i_2}-2; \cdots,p_{i_m}\;\nparallel
0,p_{i_m}-2)$.

\begin{equation}
 \label{equ:L52}
   \begin{array}{lrclll}
    m\prod_{i=1}^{i_m}\limits\left[1-\frac{1}{p_i}-\frac{\widetilde{1}}{p_i}\right]
        \left[1-\frac{1}{p_j}-\frac{\widetilde{1}}{p_j}\right]\\
    =S(m,p_1\nparallel 0; p_2\nparallel 0,p_2-2; \cdots; p_{i_m}\nparallel 0,p_{i_m}-2; p_j\nparallel 0,p_j-2)\\
    \geq S(\lceil m(1-3/p_j)\rceil,p_1\nparallel 0; p_2\nparallel 0,p_2-2; \cdots; p_{i_m}\nparallel 0,p_{i_m}-2)\\
%    +\Delta(p_i||r_{p_i})\right\} \\
    = \left\lceil
    m\left(1-\frac{3}{p_j}\right)\right\rceil
    \prod_{i=1}^{i_m}\limits\left[1-\frac{1}{p_i}-\frac{\widetilde{1}}{p_i}\right].
 \end{array}
\end{equation}
%%%%%%%%%%%%%%%%%%%%%%

Besides, suppose that for $1<r\leq i_m$,

\begin{equation}
 \label{equ:L5t}
 \begin{array}{l}
m\prod_{i=r}^{i_m}\limits\left[1-\frac{1}{p_i}-\frac{\widetilde{1}}{p_i}\right]
    \left[1-\frac{1}{p_j}-\frac{\widetilde{1}}{p_j}\right]
    = \left\lceil m\left(1-\frac{3}{p_j}\right)\right\rceil
    \prod_{i=r}^{i_m}\limits\left[1-\frac{1}{p_i}-\frac{\widetilde{1}}{p_i}\right]+t,\\
 \end{array}
\end{equation}
where $t\geq 0$. It means that the effect of the operator
$\left[1-\frac{1}{p_j}-\frac{\widetilde{1}}{p_j}\right]$ when
operating on $m$ is that dividing $Z=\{1,2,\cdots,m\}$ into two
set $X=\{1,2,\cdots,m'=\left\lceil
 m\left(1-\frac{3}{p_j}\right)\right\rceil\}$ and
 $X'=\{X'_1,X'_2,\cdots,X'_t, t=m-m'\geq 0\}$. From equation (\ref{equ:d32}), (\ref{equ:d33}), (\ref{equ:L1}),
(\ref{equ:L3}),  and (\ref{equ:L5t}), we have

\begin{displaymath}
\begin{array}{lr}
    m\prod_{i=r-1}^{i_m}\limits\left[1-\frac{1}{p_i}-\frac{\widetilde{1}}{p_i}\right]
    \left[1-\frac{1}{p_j}-\frac{\widetilde{1}}{p_j}\right]\\
    =m\prod_{i=r}^{i_m}\limits\left[1-\frac{1}{p_i}-\frac{\widetilde{1}}{p_i}\right]
    \left[1-\frac{1}{p_j}-\frac{\widetilde{1}}{p_j}\right]
    \left[1-\frac{1}{p_{r-1}}-\frac{\widetilde{1}}{p_{r-1}}\right]\\
    =\left( \left\lceil m\left(1-\frac{3}{p_j}\right)\right\rceil
    \prod_{i=r}^{i_m}\limits\left[1-\frac{1}{p_i}-\frac{\widetilde{1}}{p_i}\right]+t \right)
    \left[1-\frac{1'}{p_{r-1}}-\frac{1''}{p_{r-1}}\right]\\
    \geq \left(\left\lceil m\left(1-\frac{3}{p_j}\right)\right\rceil
    \prod_{i=r}^{i_m}\limits\left[1-\frac{1}{p_i}-\frac{\widetilde{1}}{p_i}\right]\right)
    \left[1-\frac{1'}{p_{r-1}}-\frac{1''}{p_{r-1}}\right]\\
    = \left\lceil m\left(1-\frac{3}{p_j}\right)\right\rceil
    \prod_{i=r-1}^{i_m}\limits\left[1-\frac{1}{p_i}-\frac{\widetilde{1}}{p_i}\right].\\
\end{array}
\end{displaymath}

Or
\begin{displaymath}
\begin{array}{lr}
    m\prod_{i=r}^{i_m}\limits\left[1-\frac{1}{p_i}-\frac{\widetilde{1}}{p_i}\right]
        \left[1-\frac{1}{p_j}-\frac{\widetilde{1}}{p_j}\right]
        \left[1-\frac{1}{p_{r-1}}-\frac{\widetilde{1}}{p_{r-1}}\right]\\
    =m\prod_{i=r}^{i_m}\limits\left[1-\frac{1}{p_i}-\frac{\widetilde{1}}{p_i}\right]
        \left[1-\frac{1}{p_j}-\frac{\widetilde{1}}{p_j}\right]
        -m\prod_{i=r}^{i_m}\limits\left[1-\frac{1}{p_i}-\frac{\widetilde{1}}{p_i}\right]
        \left[1-\frac{1}{p_j}-\frac{\widetilde{1}}{p_j}\right]
        \left[\frac{1}{p_{r-1}}\right]\\
   {\ \ }-m\prod_{i=r}^{i_m}\limits\left[1-\frac{1}{p_i}-\frac{\widetilde{1}}{p_i}\right]
        \left[1-\frac{1}{p_j}-\frac{\widetilde{1}}{p_j}\right]
        \left[\frac{\widetilde{1}}{p_{r-1}}\right]\\
    =\left\lceil m\left(1-\frac{3}{p_j}\right)\right\rceil
        \prod_{i=r}^{i_m}\limits\left[1-\frac{1}{p_i}-\frac{\widetilde{1}}{p_i}\right]
         +t\\
   {\ \ }-\left\lceil m\left(1-\frac{3}{p_j}\right)\right\rceil
        \prod_{i=r}^{i_m}\limits\left[1-\frac{1}{p_i}-\frac{\widetilde{1}}{p_i}\right]
        \left[\frac{1}{p_{r-1}}\right]
        -t\left[\frac{1'}{p_{r-1}}\right]\\
   {\ \ }-\left\lceil m\left(1-\frac{3}{p_j}\right)\right\rceil
        \prod_{i=r}^{i_m}\limits\left[1-\frac{1}{p_i}-\frac{\widetilde{1}}{p_i}\right]
        \left[\frac{\widetilde{1}}{p_{r-1}}\right]
        -t\left[\frac{1''}{p_{r-1}}\right]\\
    =\left\lceil m\left(1-\frac{3}{p_j}\right)\right\rceil
        \prod_{i=r}^{i_m}\limits\left[1-\frac{1}{p_i}-\frac{\widetilde{1}}{p_i}\right]
        \left[1-\frac{1}{p_{r-1}}-\frac{\widetilde{1}}{p_{r-1}}\right]
        +t\left[1-\frac{1'}{p_{r-1}}-\frac{1''}{p_{r-1}}\right]\\
    \geq \left\lceil m\left(1-\frac{3}{p_j}\right)\right\rceil
        \prod_{i=r-1}^{i_m}\limits\left[1-\frac{1}{p_i}-\frac{\widetilde{1}}{p_i}\right].
\end{array}
\end{displaymath}

In fact, if
\begin{displaymath}
   \begin{array}{lr}
    m\prod_{i=1}^{i_m}\limits\left[1-\frac{1}{p_i}-\frac{\widetilde{1}}{p_i}\right]
    \left[1-\frac{1}{p_j}-\frac{\widetilde{1}}{p_j}\right]
    < \left\lceil m\left(1-\frac{3}{p_j}\right)\right\rceil
    \prod_{i=1}^{i_m}\left[1-\frac{1}{p_i}-\frac{\widetilde{1}}{p_i}\right],
\end{array}
\end{displaymath}
then for any $p_i, p_j$, before deleting the multiples of other primes, it must have \\
\begin{displaymath}
   \begin{array}{lr}
    m\left[1-\frac{1}{p_i}-\frac{\widetilde{1}}{p_i}\right]
    \left[1-\frac{1}{p_j}-\frac{\widetilde{1}}{p_j}\right]
    < \left\lceil m\left(1-\frac{3}{p_j}\right)\right\rceil
    \left[1-\frac{1}{p_i}-\frac{\widetilde{1}}{p_i}\right].
\end{array}
\end{displaymath}
which contradicts  Lemma \ref{lem:3}. So this lemma is true.

If $p_{r-1}=2$, then
\begin{equation}
 \label{equ:L53}
   \begin{array}{lr}
    m\left[1-\frac{1}{p_{r-1}}\right]\prod_{i=r}^{i_m}\limits\left[1-\frac{1}{p_i}-\frac{\widetilde{1}}{p_i}\right]
        \left[1-\frac{1}{p_j}-\frac{\widetilde{1}}{p_j}\right]\\
    \geq \left\lceil m\left(1-\frac{3}{p_j}\right)\right\rceil\prod_{j}
        \left[1-\frac{1}{p_{r-1}}\right]
        \prod_{i=r}^{i_m}\limits\left[1-\frac{1}{p_i}-\frac{\widetilde{1}}{p_i}\right].\\
\end{array}
\end{equation}
\end{proof}

With Lemma\ \ref{lem:3}, the operator
$\left[1-\frac{1}{p_j}-\frac{\widetilde{1}}{p_j}\right]$ can be
represented by $\left(1-\frac{3}{p_j}\right)$, and other operator
$\left[1-\frac{1}{p_i}-\frac{\widetilde{1}}{p_i}\right]$ can
operate on this inequality unchanged.

\section{Explantation}
\label{sect:Explantation}

Let $m=sp_ip_j+ap_j+b\geq p_j^2$, then for all $p_i < p_j$,
$m\left[1-\frac{1}{p_j}-\frac{\widetilde{1}}{p_j}\right]$ is
effective to a nature sequence $X=\left\{f(Z)\right\}$ whose
number is not less than $\lceil m(1-3/p_j)\rceil$. The reason is
as follows. When a nature sequence is deleted by the items of $Z_k
\mod p_j=0,p_j-2$, the sequence is subtracted by
$\left[\frac{m}{p_j}\right]+\left[\frac{m+2}{p_j}\right]$. We can
arrange the m items in a table of $p_j$ rows (Table~1).
$\left[\frac{m}{p_j}\right]$ will delete the $p_j th$ row, and
$\left[\frac{m+2}{p_j}\right]$ will delete the $(p_j-2)th$ row.
Thus there are $(p_j-2)$ rows left in which each item $Z_k \mod
p_j\ne0, p_j-2$.

\begin{table}[ht]
 \label{table1} \caption{Set Z}
\renewcommand\arraystretch{1.5}
 \noindent\[
%\begin{array}{|c|c|@{\cdots}|c||@{\cdots}|c||c|@{\cdots}|c|}
\begin{array}{|cc@{\ \cdots\ }c|@{\ \cdots\ }c|c@{\ \cdots\ }c|}
\hline
  1 &  p_j+1 &  (p_i-1)p_j+1 &  (sp_i-1)p_j+1 & sp_ip_j+1  & sp_ip_j+ap_j+1\\
%\hline
  2 &  p_j+2 &  (p_i-1)p_j+2 &  (sp_i-1)p_j+2 & sp_ip_j+2  & sp_ip_j+ap_j+2\\
%\hline
  \cdots &  \cdots &  \cdots &  \cdots & \cdots  & \cdots\\
%\hline
  b &  p_j+b &  (p_i-1)p_j+b &  (sp_i-1)p_j+b & sp_ip_j+b  & sp_ip_j+ap_j+b\\
%\hline
  \cdots &  \cdots &  \cdots &  \cdots & \cdots &\\
%\hline
  p_j &  2p_j &  p_ip_j &  sp_ip_j & sp_ip_j+p_j &\\
\hline
\end{array}
\]
\end{table}

But every $p_i$ items ($0\leq Z_k \mod p_i\leq p_i-1$) in any row
of the first $sp_i$ columns consist in a complete system of
residues modulo $p_i$, because $C_1=\{1, p_j+1, 2p_j+1, \cdots ,
(p_i-1)p_j+1\}$ and $C_r=\{C_1+r\}$  are both complete system of
residues modulo $p_i$, where $r$ is any (row or  column) constant.
There are $(p_j-2)$ such rows or $sp_i(p_j-2)$ items left. These
items are effective to a nature sequence when deleting the items
of $Z_k \mod p_i=0,p_i-2$.

$sp_ip_j\left[1-\frac{1}{p_i}-\frac{\widetilde{1}}{p_i}\right]
   \left[1-\frac{1}{p_j}-\frac{\widetilde{1}}{p_j}\right]
   =sp_i(p_j-2)\left[1-\frac{1}{p_i}-\frac{\widetilde{1}}{p_i}\right]=s(p_i-2)(p_j-2).$

Let $t=ap_j+b$, $0 \leq b \leq p_j-1$, $0 \leq a \leq p_i-1$.
$t\left[\frac{1}{p_j}+\frac{\widetilde{1}}{p_j}\right]$ will
delete at most $a+(a+1) \leq t\frac{1}{p_j}+sp_ip_j\frac{1}{p_j}=
\frac{m}{p_j}$ items. If we add these items by removing those from
the end of sequence then the sequence is again effective to a
nature sequence, the sequence left has at least,
\begin{displaymath}
\begin{array}{rl}
    M(j) & \geq sp_i(p_j-2)+t\left[1-\frac{1}{p_j}-\frac{\widetilde{1}}{p_j}\right]-[a+(a+1)]\\
     & \geq sp_i(p_j-3)+t+sp_i-4a-2 \\
     & =sp_i(p_j-3)-3\left\lfloor\frac{t}{p_j}\right\rfloor+t+(sp_i-a-2)\\
     & \geq sp_ip_j(1-\frac{3}{p_j})+\left\lceil t-\frac{3t}{p_j}\right\rceil+(sp_i-a-2)\\
     & = \left\lceil(sp_ip_j+t)(1-\frac{3}{p_j})\right\rceil+(sp_i-a-2)\\
     & = \left\lceil m(1-\frac{3}{p_j})\right\rceil+(sp_i-a-2).
\end{array}
\end{displaymath}

For $s \geq 2$ or $a\leq p_i-2$, we have
$(sp_i-a-2)=(s-1)p_i+(p_i-a-1)-1\geq 0$, $M(j) \geq \left\lceil
m(1-\frac{3}{p_j})\right\rceil$.

For $s =1$ and $a=p_i-1$, the items of t have $p_j$ rows, $p_i-1$
columns and some $b$ items. In each of the first $b$ rows, there
are exact $p_i$ items which consist in a complete system of
residues modulo $p_i$, and these items can be considered as an
effective nature sequence when deleting the multiples of $p_i$
($Z_k \mod p_i=0,p_i-2$). The other items have at most $p_j$ rows
and $p_i-1$ columns where the multiples of $p_j$ have at most
$2(p_i-1)$. As before, we can add these items to make the $t$ as
an effective nature sequence, therefore,
\begin{displaymath}
\begin{array}{rl}
    M(j) &\geq
    sp_i(p_j-2)+t\left[1-\frac{1}{p_j}-\frac{\widetilde{1}}{p_j}\right]-2(p_i-1)\\
    & \geq  \left\lceil m(1-\frac{3}{p_j})\right\rceil+(sp_i-a-1)\geq
    \left\lceil m(1-\frac{3}{p_j})\right\rceil.
\end{array}
\end{displaymath}
Thus for any $p_i<p_j$, the original sequence of $m\geq p_j^2$,
when deleted by the items of $Z_k \mod p_j=0,p_j-2$ from $Z$, is
effective to reconstruct a new natural sequence having at least
$\left\lceil m(1-\frac{3}{p_j})\right\rceil$ items.

\begin{example} $n=41$, $P=\{2,3,5\}$, $Z=\{1,2,\cdots,41\}$.

For $p_j=5$, after deleted the item of $Z_k \mod p_j=0,5-2$, it becomes\\
$Z\rightarrow
Z'=\{1,2,4,6,7,9,11,12,14,16,17,19,21,22,24,26,27,29,31,32,34,36,\-37,\-39,41\}$.
we can rearrange these items as
$Z'=\{1,2,(21),4,(29),6,7,(26),9,\-(22),\-11,12,(19),14,(27),16,17||,24,31,32,34,36,\-37,\-39,41\}$.
The first $\left\lceil
n\left(1-\frac{3}{p_j}\right)\right\rceil=17$ items can be taken
as an effective natural sequence from the original one when
deleting the items of $Z_k\mod3=0,3-2$. The other sequence
$X=\{24,31,32,34,36,\-37,\-39,41\}$, having at least zero item
when deleting the  multiples of all primes, will be neglected in
further procession.
\end{example}

%%%%%%%%%%%%%%%%%%%%%%%%%%%%%%%%%5

\section{ Proof of Theorem~\ref{theorem}}
\label{sect:Theorem}

For a given $n$, consider the possible twin primes of $[Z_k,
Z_k+2]$ in set $Z=\{1,2,\cdots,n\}$, $p_v^2\leq n<p_{v+1}^2$.

\begin{lemma}
\label{thm:lem2} The number of twin prime pairs in $n$,
\begin{equation}
  \label{equ:TL2}
   \begin{array}{rl}
    D(n) \geq\left\lceil\frac{p_v}{3}\prod_{i=3}^{v-1}\frac{p_{i+1}-3}{p_i}\right\rceil
            +D(\sqrt{n})-2.
 \end{array}
\end{equation}
\end{lemma}
\begin{proof}
Because $m(v)=n\geq p_v^2$. Let $m(j-1)=\left\lceil
m(j)(1-3/p_j)\right\rceil$, then
\begin{displaymath}
   \begin{array}{rl}
   m(j-1) & =\left\lceil m(j)(1-3/p_j)\right\rceil
         \geq p_j^2(1-3/p_j)=p_j(p_j-3)     \\
        &  \geq (p_{j-1}+2)(p_{j-1}-1)=p_{j-1}^2+p_{i}-2
          \geq p_{j-1}^2.
   \end{array}
\end{displaymath}
And for any $i\leq (j-1)$, we have $m(i)\geq p_i^2$.

From equation (\ref{equ:D0n}), (\ref{equ:L5}), (\ref{equ:L53}) and
(\ref{equ:L4}),

\begin{displaymath}
   \begin{array}{llr}
    D_0(n) & =n\left[1-\frac{1}{2}\right]\prod_{i=2}^v\limits\left[1-\frac{1}{p_i}-\frac{\widetilde{1}}{p_i}\right]\\
        & \geq \left\lceil n\left(1-\frac{3}{p_v}\right)\right\rceil\left[1-\frac{1}{2}\right]
         \prod_{i=2}^{v-1}\limits\left[1-\frac{1}{p_i}-\frac{\widetilde{1}}{p_i}\right]\\
        & \geq \left\lceil\left\lceil n\left(1-\frac{3}{p_v}\right)\right\rceil
         \left(1-\frac{3}{p_{v-1}}\right)\right\rceil
         \left[1-\frac{1}{2}\right]
        \prod_{i=2}^{v-2}\limits\left[1-\frac{1}{p_i}-\frac{\widetilde{1}}{p_i}\right]\\
       & \geq \cdots\\
     & \geq \left\lceil n\prod_{i=3}^{v}\limits\left(1-\frac{3}{p_i}\right)\right\rceil
    \left[1-\frac{1}{2}\right]
        \left[1-\frac{1}{p_2}-\frac{\widetilde{1}}{p_2}\right]\\
     & \geq \left\lceil\frac{n}{6}\prod\limits_{i=3}^{v}\left(1-\frac{3}{p_i}\right)\right\rceil-1\\
 \end{array}
\end{displaymath}

but
\begin{displaymath}
   \begin{array}{lr}
   \quad \prod\limits_{i=3}^v\left(1-\frac{3}{p_{i}}\right)
   \ =\frac{p_3-3}{p_3}\frac{p_4-3}{p_4}\cdots\frac{p_{v-1}-3}{p_{v-1}}\frac{p_v-3}{p_v}\\
    \quad\ = \frac{p_3-3}{p_v}\frac{p_4-3}{p_3}\frac{p_5-3}{p_4}\cdots\frac{p_{v}-3}{p_{v-1}}
   \ = \frac{2}{p_v}\prod\limits_{i=3}^{v-1}\frac{p_{i+1}-3}{p_i},
\end{array}
\end{displaymath}

and $n\geq p_v^2$, so
\begin{displaymath}
   \begin{array}{l}
    \quad
    D_0(n)\geq\left\lceil\frac{n}{3p_v}\prod_{i=3}^{v}\left(1-\frac{3}{p_i}\right)\right\rceil-1
       \geq\left\lceil\frac{p_v}{3}\prod_{i=3}^{v-1}\frac{p_{i+1}-3}{p_i}\right\rceil-1.
\end{array}
\end{displaymath}

From equation (\ref{equ:Dn}),(\ref{equ:D1}),

\begin{displaymath}
   \begin{array}{rl}
    D(n) & =D_0(n)+D(\sqrt{n})-D_1\\
          & \geq\left\lceil\frac{p_v}{3}\prod_{i=3}^{v-1}\frac{p_{i+1}-3}{p_i}\right\rceil-1+D(\sqrt{n})-D_1\\
          & \geq\left\lceil\frac{p_v}{3}\prod_{i=3}^{v-1}\frac{p_{i+1}-3}{p_i}\right\rceil
            +D(\sqrt{n})-2.
 \end{array}
\end{displaymath}
\end{proof}

\begin{proof}[{\bf Proof of Theorem \ref{theorem}}]
Suppose that there is no twin prime when greater than enough large
number $n_M$.

\begin{equation}
  \label{equ:TP1}
   \begin{array}{rl}
    D(n)  & \leq D(n_M) \quad \mbox{for enough}\ n>n_M.
 \end{array}
\end{equation}

Consider the pairs of twin prime in the range $[1, n]$, where
$n=n_M^2$,

\begin{equation}
  \label{equ:TP2}
   \begin{array}{rl}
    D(n)  & \geq\left\lceil\frac{p_V}{3}\prod_{i=3}^{V-1}\frac{p_{i+1}-3}{p_i}\right\rceil
            +D(\sqrt{n})-2.
 \end{array}
\end{equation}
where $p_V$ is the maximum prime in $n$.

Then the twin prime between $n_M$ and $n=n_M^2$ will have,
\begin{equation}
  \label{equ:TP3}
   \begin{array}{rl}
    \Delta D(n) & =D(n)-D(n_M)=D(n)-D(\sqrt{n})\\
                & \geq\left\lceil\frac{p_V}{3}\prod_{i=3}^{V-1}\frac{p_{i+1}-3}{p_i}\right\rceil-2
                 =\left\lceil \frac{W(V)}{3}\right\rceil-2.
 \end{array}
\end{equation}
where
\begin{equation}
  \label{equ:TP4}
    W(V) = p_V\prod_{i=3}^{V-1}\frac{p_{i+1}-3}{p_i}.
\end{equation}

If $W(V)>6$ then $\Delta D(n)\geq 1$, there will be at least one
twin prime between $n_M$ and $n=n_M^2$.

For $V\geq V_0=5, p_{V_0}=11$,

\begin{equation}
  \label{equ:TP5}
    W(5)= 11\frac{4}{5}\frac{8}{7}=10.057>6
\end{equation}

Suppose that for $V$,
$W(V)=p_V\prod_{i=3}^{V-1}\frac{p_{i+1}-3}{p_i}> 6$, then for
$V+1$,

\begin{equation}
  \label{equ:TP6}
   \begin{array}{rlr}
   W(V+1) & =p_{V+1}\prod_{i=3}^{V-1}\frac{p_{i+1}-3}{p_i}\frac{p_{V+1}-3}{p_V}\\
          & =p_{V}\prod_{i=3}^{V-1}\frac{p_{i+1}-3}{p_i}\frac{p_{V+1}}{p_V}\frac{p_{V+1}-3}{p_V}\\
          & =W(V)\frac{p_{V+1}^2-3p_{V+1}}{p_V^2}.\\
\end{array}
\end{equation}
Because, $p_{V+1}=p_V+2\Delta, \Delta\geq 1, p_V\geq 5$,

\begin{displaymath}
   \begin{array}{lr}
   p_{V+1}^2-3p_{V+1}-p_V^2
    =(p_V+2\Delta)^2-3(p_V+2\Delta)-p_V^2\\
    =p_V^2+4\Delta p_V+4\Delta^2-3p_V-6\Delta-p_V^2\\
    =3(\Delta-1) p_V+\Delta (p_V+4\Delta-6)
    > 0.
\end{array}
\end{displaymath}

Therefore, $\frac{p_{V+1}^2-3p_{V+1}}{p_V^2}>1$, and
\begin{equation}
  \label{equ:TP7}
   \begin{array}{lr}
    W(V+1)>W(V)>6.\\
\end{array}
\end{equation}
So that $\Delta D(n)\geq 1$. It contracts the supposition of Eq.
(\ref{equ:TP1}). Therefore `there are infinitely many pairs of
twin prime'. From Eq. (\ref{equ:TP3}) and (\ref{equ:TP4}), $\Delta
D(n)$ approaches infinity as n grows without bound. The proof is
completed.
\end{proof}

\begin{example}[Actual vs. Simplified Formula] \label{sect:figure}

Let
\begin{equation}
  \label{equ:F1}
   \begin{array}{rl}
    D'(n)  & =\left\lceil\frac{p_v}{3}\prod_{i=3}^{v-1}\frac{p_{i+1}-3}{p_i}\right\rceil-2.\\
 \end{array}
\end{equation}
Then from Eq. (\ref{equ:TP2}), it should have $D(n)\geq D'(n)$.

Figure ~\ref{Fig:Gfig1} shows the actual pairs of twin prime
$D(n)$ (solid line) in the range of $[1, p_v^2+1]$, and its
simplified formula $D'(n) $(dashed line) from Eq. (\ref{equ:F1}).
It clearly shows that $D(n)\geq D'(n)$, and $D'(n)$ has no up
bound for enough large $n$.
%------------------------------- Fig1: MinTwins
\begin{figure}[ht]
\begin{center}
    \includegraphics[scale=0.4,width=120mm,height=60mm,,angle=0]{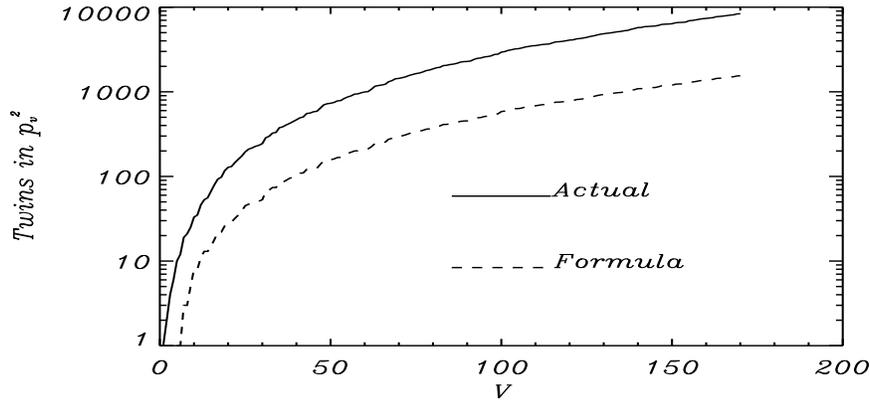}
\end{center}
\caption{The minimum number of actual prime pairs (solid) and its
simplified  formula (dashed) against v.}
   \label{Fig:Gfig1}
\end{figure}

%------------------------------- Fig2: ratio

\begin{figure}[ht]
\begin{center}
    \includegraphics[scale=0.4,width=120mm,height=60mm,,angle=0]{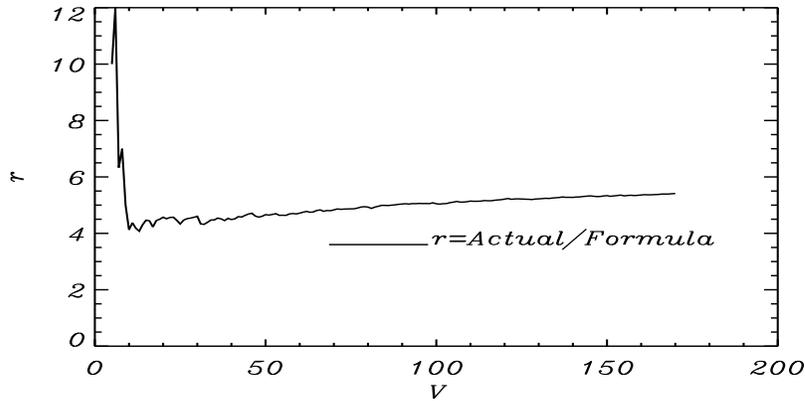}
\end{center}
\caption{The ratio of $r=D(n)/D'(n)$ against v.}
   \label{Fig:Gfig2}
\end{figure}
Figure ~\ref{Fig:Gfig2} shows the ratio of
\begin{equation}
  \label{equ:F2}
   \begin{array}{rl}
    r=\frac{D(n)}{D'(n)}.
 \end{array}
\end{equation}
It shows that, this ratio is greater than one and seems to
progress as $n$ increase. It also shows that the formula (Eq.
\ref{equ:F1}) used in our proof is only a small part of the actual
twin prime pairs.

\end{example}

\bibliographystyle{99}%{amsplain}

\end{document}